\documentclass[leqno]{article}
\def\cyr{\tencyr\cyr}              
\usepackage[all]{xy}
\usepackage{amsmath,amscd,amssymb,amsthm,amsxtra,latexsym,epsfig,epic,eepic,mathrsfs,makeidx}

\DeclareMathOperator{\es}{es}
\DeclareMathOperator{\ES}{ES}
\DeclareMathOperator{\eb}{eb}

\DeclareMathOperator{\IN }{IN}

\DeclareMathOperator{\Hom}{Hom}

\DeclareMathOperator{\ord}{ord}

\DeclareMathOperator{\Spec}{Spec}

\DeclareMathOperator{\SR}{SR}

\DeclareMathOperator{\AES}{AES}

\newcommand{\N}{{\mathbb N}}

\newcommand{\C}{{\mathbb C\,}}

\renewcommand{\tilde}{\widetilde}

\newcommand{\ko}{{\mathcal O}}

\newcommand{\lra}{\longrightarrow}

\newcommand{\tensor}{\otimes}
\newcommand{\eps}{\varepsilon}

\newcommand{\vfi}{\varphi}

\renewcommand{\bar}{\overline}

\newcounter{pren}
\newtheorem{Lemma}{Lemma}[section]
\newtheorem{Proposition}[Lemma]{Proposition}
\newtheorem{Theorem}[Lemma]{Theorem}

\theoremstyle{definition}
\newtheorem{defn}[Lemma]{Definition}
\newtheorem{example}[Lemma]{Example}

\newtheorem{remark}[Lemma]{Remark}

\begin{document}
\title{Equisingular Deformations
of Plane Curve and of Sandwiched Singularities}
\date{}
\author{Theo de Jong\\
email: dejong@math.uni-sb.de\\}
\maketitle

\begin{abstract}
Let $C$ be an isolated plane curve singularity, 
and $(C,l)$ be a decorated curve. In 
this article we compare the equisingular
deformations of $C$ and the sandwiched
singularity $X(C,l)$. We will  prove that
for $l \gg 0$ the functor of equisingular deformations
of $C$ and $(C,l)$ are equivalent. 
From this we deduce a proof of a formula for the dimension
of the equisingular stratum.
 Furthermore we will 
show how compute  the equisingularity ideal
of the curve singularity $C$,  given the minimal (good) resolution
 of $C$. 
\end{abstract}

\section{Introduction}

In the sixties, Zariski \cite{z1}, \cite{z2}, \cite{z3} started the modern 
study of equisingularity of plane curve singularities.
Two plane curve singularities are called equisingular if
one can simultaneously resolve their singularities.
Zariski proved that equisingular plane curve singularities
are topologically equivalent. In particular, the characteristic
exponents of the branches are the same, 
and corresponding intersection numbers are the same.
For more characterizations, we refer
to the article of Teissier, \cite{tei}.

Infinitesimal equisingular families of plane curve singularities
were studied by Wahl, \cite{wahl}. He defined 
the functor of Artin rings  of equisingular 
deformations of plane curve singularities,
and showed that it is a smooth subfunctor
of deformations of the plane curve singularity itself. 
He showed that the tangent space can be 
identified with an ideal in the local ring of the curve:
the so-called equisingularity ideal. 
For an irreducible curve singularity, it 
seems, given a parametrization, not to be to
difficult to compute  this equisingularity ideal.
see \cite{tei}.

In the fourth section of this article we will give an (inductive) algorithm for
computing the equisingularity ideal, given  
a minimal good resolution of the curve singularity.
The idea is to look at a {\em simultaneous resolution} functor.
This is a functor (of Artin rings) which describes 
all deformations of the plane curve singularity (with sections) which
can be simultaneously resolved. 
If one considers such deformations over smooth
curves, the singularity
is allowed to split up in several singularities.
For the equisingular deformation functor, this is
not allowed.  
The tangent space of this simultaneous resolution
functor can be calculated inductively. 
An infinitesimal deformation with simultaneous resolution
in particular gives an infinitesimal deformation of the 
minimal embedded good resolution of $C$. In particular we get for each 
exceptional curve a deformation of a small neighborhood of
this curve. If all these deformations are 
trivial we have that are deformation is equisingular,
by using results of Wahl \cite{wahl}.
Thus the equisingularity ideal is
the kernel of a map $\vfi$ from the tangent
space of the functor of deformations with simultaneous
resolution to a direct sum ofcohomology groups $\oplus H^1(F_i,N_{F_i/Z})$.
Here $Z$ is the minimal good resolution, and the
$F_i$ are the exceptional curves of $Z$.

The study of equisingular deformations of {\em surface
singularities} was started by Wahl in \cite{wahl1}.
This theory is more difficult than for curve singularities.
For rational surface singularities, however,
one still has a good theory. In this
paper we study equisingular deformations
of a  special, but nevertheless rather broad,
 class of rational surface
singularities: the sandwiched singularities. These
surface singularities were studied by
Zariski \cite{z4}, Lipman \cite{lip}, Hironaka \cite{hir},
Spivakovsky \cite{spiva}, De Jong and Van Straten
\cite{djvs2} and Gustavsen \cite{G}.
These sandwiched singularities can be constructed by 
so--called decorated curves $(C,l)$ see \cite{djvs2},
so are of type $X(C,l)$. 
To each branch of $C$ one assigns a natural number 
satisfying a technical condition. 
We first show that the analytic type of $C$ is
uniquely determined by the analytic type of $X(C,l)$,
as soon as $l \gg 0$. Then 
we will show that for $l\gg 0$ the equisingular
deformations of $C$ are in one-one correspondence
with equisingular deformations of the sandwiched
singularity $X(C,l)$. 
Using two different formulas
for the dimension of  the Artin component of $X(C,l)$,
we are then able to deduce a formula for the
dimension of the equisingular stratum of $C$.

\section{Equisingular Deformations of  Sandwiched
Singularities.}

We will  consider {\em sandwiched} singularities.
First consider a curve singularity $C = \cup_{i \in T}C_i$.
We let $m(i)$ be the sum of the multiplicities
of branch $i$ in the multiplicity sequence of the minimal embedded
resolution of $C_i$, and $M(i)$ be the sum of the multiplicities
of branch $i$ in the multiplicity sequence of the minimal 
good embedded
resolution of $C_i$. Let  $\eb(C)$ be  the number of extra blowing-ups needed
to come from a {\em minimal embedded} resolution of $C$ to
a {\em minimal {\bf good} embedded} resolution of $C$. 
 Note the following: 
\begin{Lemma}\label{propa}
\[
\sum_{i \in T} M(i) -m(i) = \eb(C).
\]
\end{Lemma}
Let $l: T \lra \N$ be such that $l(i) \ge M(i) +1$ for all $i$. 
The pair $(C,l)$ is called a {\em decorated curve}, 
see \cite{djvs2}, 1.3. 
One gets the {\em sandwiched singularity}
$X(C,l)$ as follows. First take a minimal good embedded resolution
of $C$, and then do $l(i)-M(i)$ consecutive blowing-ups
at the strict transform of the $i$'th branch, thereby inducing
a chain  of $l(i) -M(i)$ new exceptional curves for each $i$. 
 We get a modification
of $\C^2$:
\[
(Z(C,l),F)\lra (\C^2,0).
\]
Consider $E$,  the subgraph of all irreducible components whose
self-intersection is not $-1$. From the resolution process 
one sees that $E$ is connected, has negative definite intersection
matrix, and so can be contracted by the result of Grauert--Mumford
to a normal surface singularity $X(C,l)$, which one easily 
sees to be rational. This singularity is denoted by $X(C,l)$.
The modification $(Z(C,l),F)$ we can get by blowing up
a {\em complete ideal} $I(C,l)$.\\

The following theorem
 was proved by Gustavsen \cite{G} for 
the case the
plane curve singularities are irreducible.

\begin{Theorem}\label{th27900}
Consider decorated curves $(C,l)$ and $(C',l')$. Suppose that
the surface singularities 
 $X(C,l)$ and $X(C',l')$ are isomorphic.
Suppose $l \gg 0$ and $l' \gg 0$.
Then $C$ and $C'$ are isomorphic and $l = l'$.
(This is to be interpreted that for some isomorphism 
corresponding branches have the same number attached.)
\end{Theorem}

\begin{proof}
Let $\vfi$ be an isomorphism taking $X(C',l')$
to $X(C,l)$. By the universal property of blowing--up,
this extends to an isomorphism between the minimal 
resolutions $\vfi: (\tilde{X}(C',l'), E') \lra (\tilde{X}(C,l),E)$.
In particular the dual graphs of $X(C,l)$ and $X(C',l')$
are isomorphic.
We write $C = \cup_{i \in T} C_i$.
As both $l\gg 0 $ and $l'\gg 0$
there are exactly $\#T$ very long chains 
of $(-2)$--curves in both $E$ and $E'$. 
Due to their length, they can be recognized. To each endpoint 
$G'_i$ of such a very long chain in $E'$
there is, by the  construction
of sandwiched singularities, a $(-1)$--curve $F_i'$ in $Z(C',l')$ and a smooth 
curve $C_i'$ intersecting $F'_i$. Similar for $Z(C,l)$.
The isomorphism $\vfi$ sends  $G'_i$ to a $G_i$,
as the resolution graphs of $X(C',l')$ and $X(C,l)$ are isomorphic.  
By using the isomorphism $\vfi$, we may glue every $F_i'$
(and in it $C_i'$) to the curve $G_i$ in $X(C,l)$. 
Thus we may assume that the 
two resolutions are equal: $(Z(C',l'), E') = (Z(C,l),E)$.

Without loss of generality, we may assume that
both $x$ and $y$ are generic elements of $\ko_{C'}$.
By looking at the image under $\vfi$ of the divisors of the pull-back
of $(x)$ and $(y)$ on $X(C',l')$ (which is the fundamental cycle of
$Z(C'l')$ restricted to $E'$ plus a non-compact curve
intersecting the first blown up curve), we
see that the divisors of  $\vfi(x)$ and $\vfi(y)$ are the 
divisors of two functions which generate the maximal ideal
of $\ko_C$. Thus $\vfi$  maps $C'$ to an isomorphic curve,
which except for the $(-1)$--curves has the same
resolution as $C$. Thus  we may and will suppose that
the  resolutions  $Z(C',l')$ and $ Z(C,l)$ 
are equal.
From the construction of the $Z(C,l)$ out of
the $\tilde{X}(C,l)$ it 
follows that $l = l'$.

There is an algorithm to get, from the resolution
of a plane curve singularity, the resolution
of every irreducible component.
Following this algorithm,
we see that in the resolution of any irreducible component
 $C_i$ of $C$ there is a very long
chain of $(-2)$--curves. By the remarks above
this is, except for the $(-1)$--curves, also 
the resolution for $C_i'$.
By the formula for the intersection number,
see, for example, \cite{JP} Theorem 5.4.8,
the intersection number between $C$ and $C'$ increases
if the chains of $(-2)$ curves, and thus $l(i)$ and $l'(i)$
become bigger. This intersection number is equal
to both  $\dim \C\{x,y\}/(f_i,f_i')$, 
and also to the vanishing order of $f_i'$ on $\tilde{C_i}$,
where $\tilde{C_i}$ is the normalization of $C_i$.
Here the curves $C_i$ and $C_i'$ are defined
by irreducible $f_i$ and $f_i'$.  
Thus we may assume that for all $i$ the vanishing order of 
$f_i'$ on $\tilde{C_i}$ is at least $k\cdot c(i)$, for 
some large $k$, and where $c(i)$ is the $i$'th conductor
number of $C$, corresponding to the branch $C_i$. 
By definition of the conductor, every function which
vanishes with order at least $c(i)$ on the normalization
$\tilde{C}_i$ for all $i$,  is an element
of the maximal ideal $(x,y)$ of $\ko_C = \C\{x,y\}/(f)$.
It follows that the class of $f' = \prod f_i'$ in $\ko_C$
lies in $(x,y)^k$. Thus $uf - f' \in (x,y)^k$,
where we now view $(x,y) $ as an ideal in $\C\{x,y\}$. 
By symmetry, $u$ is a unit. 
By taking $k$ large, we may, by the finite determinacy 
theorem, assume that $uf$ and $f'$ are right equivalent.
In particular, their zero sets $C$ and $C'$ 
are isomorphic. This is what we had to show. 
\end{proof}

We now study equisingular deformations.
So we have an functor $\ES_C$ of equisingular deformation of
the curve singularity $C$, and the functor 
$ES_{X(C,l)}$ of equisingular deformations of
the sandwiched singularity $X(C,l)$. 
The dimension of the Zariski tangent spaces we denote by 
$\es(C)$ and $\es(X(C,l))$. 
We quote the following result due
to Gustavsen, see \cite{G}, Theorem 3.3.22.

\begin{Theorem}
There is a natural formally smooth map of functors
\[
\ES_C \lra \ES_{X(C,l)}.
\]
\end{Theorem}

This in particular, for any equisingular deformation
of $X(C,l)$ we can find an equisingular deformation of $C$,
mapping to it naturally, but this equisingular deformation of $C$
might not be unique, even on tangent spaces.
The existence of the map $\ES_C \lra \ES_{X(C,l)}$ is quite obvious.
An equisingular deformation 
of $C$ induces a deformation of $Z(C,l)$,
and thus a deformation of the resolution of $X(C,l)$. This
thus gives an equisingular deformation of $X(C,l)$. 

Our aim is to prove that for $l\gg 0$, 
the formally smooth map is in fact an isomorphism.
For this, it suffices to show that it is an isomorphism
on tangent spaces.

\begin{Proposition}\label{Prop27900}
Let $C$ be an isolated curve singularity, and $(C,l)$ 
be a decorated curve. Suppose that $l \gg 0$. 
Then the Zariski tangent spaces of the equisingular
deformations of $X(C,l)$ and of $C$ are isomorphic. In particular
for the dimensions we have
$\es(X(C,l)) = \es(C).$
\end{Proposition}

We use the theory 
developed in \cite{djvs2}. 
Let $C$ be given by $f = 0$, and $c(i)$ be the conductor
number on branch $C_i$. We consider a function
$g(x,y)$ whose vanishing order on the normalization
of $C_i$ is equal to $c(i) + l(i)$ for all $i \in T$. 
In \cite{djvs2} the following non--isolated surface singularity
in  $\C^3$ is  considered:
\[
Y(C,l) = V(zf -g).
\]
on the choice of $g$. 
The normalization of $Y(C,l)$ is proved
to be $X(C,l)$. 
Let $\Sigma \subset \C^2$ be the space defined by the
conductor $I$ of $C$. We will consider R.C.--deformations
of the pair $(\Sigma,C)$. For the definition of R.C.--deformations,
see for example
\cite{djvs2}, Appendix,  and the references mentioned there.
This functor is canonically isomorphic to 
deformations of the diagram  $\tilde{C} \lra C$,
where   $\tilde{C} \lra C$ is the normalization. 
The space of infinitesimal deformation 
of R.C.--deformations of the  pair $(\Sigma,C)$,
we denote by $T^1(\Sigma,C)$.
It was stated, but not properly shown in \cite{djvs2}, Remark 3.18,
that for $l \gg 0$ we have an exact sequence
\begin{equation}\label{eq28900}
0 \lra I^{ev}/(f,\Theta_\Sigma(g)) \lra T^1_{X(C,l)} \lra T^1(\Sigma,C)
\lra 0.
\end{equation}
Here $I^{ev} = \{g \in \C\{x,y\}:   \ord(g_{| \tilde{C}_i})\ge l(i) +
  m(i)\; {\rm for \; all\;} i\}$ 
 for $\tilde{C}_i \lra C_i$ the normalization. 
Furthermore $\Theta_\Sigma$ are all derivations $\theta$
of $\C\{x,y\}$ satisfying the two conditions $\theta(f) \subset (f)$ and
$\theta(I) \subset I$. If $t_i$ is a parameter on $\tilde{C}_i$,
one can identify $\Theta_\Sigma$ with 
the module generated by the $t_i \frac{\partial}{\partial
  t_i}$. 
Let us see how the sequence \ref{eq28900} comes about. In \cite{djvs2},
Proposition 3.7 it is shown that every infinitesimal deformation
of $X(C,l)$
can be obtained, up to isomorphism, as follows. Take an 
infinitesimal R.C.--deformation 
$(\Sigma_\eps,C_\eps)$ of  $(\Sigma, C)$, and a $g_\eps$ such
that also $(\Sigma_\eps, (g_\eps =0))$ is an R.C--deformation
of $(\Sigma, (g =0))$. Suppose that $\Sigma_\eps$ is defined
by the ideal $I_\eps$, and $C_\eps$ by $f_\eps = 0$. We get a space
$Y_\eps$ defined by $zf_\eps - g_\eps = 0$, 
and the space $X_\eps$  is defined by the local
ring 
\[
\Hom_{Y_\eps}(I_\eps,I_\eps).
\]
So $X_\eps$ is obtained by ``simultaneous normalization''.
In general, the choice of the isomorphism class of the
R.C.--deformation of $(\Sigma,C)$ is not uniquely determined by the
given infinitesimal deformation of $X(C,l)$. But we 
will see  and have to show that in case $l \gg 0$ it 
in fact does. By linearity, it suffices to 
show that if we have a trivial infinitesimal deformation
of $X(C,l)$, we are only allowed to deform $(\Sigma, C)$
trivially. Now a trivial first order family can be 
extended to a trivial family over a germ of a smooth curve $T$, 
and we may assume that this family is a product family. 
By the results of \cite{djvs2}, in particular Theorem 3.3
we can therefore extend
the R.C--deformation $(\Sigma_\eps, C_\eps)$ and 
$(\Sigma_\eps,(g_\eps = 0))$  to  deformations over this smooth curve. 
Thus we get an $I_T$ defining $\Sigma_T$, $f_T$ defining
$C_T$,  and some $g_T$, such that the local ring of the  trivial 
family is given
by $\Hom_{Y_T}(I_T,I_T)$,
where $Y_T$ is given by $zf_T-g_T =0$.

Writing $F_t = \sum_it^if_i$, we have that for all
small $t$ the zero set of $F_t$ intersects on $Z(C,l)$ the 
$(-1)$--curve. This is because the induced deformation
of $X(C,l)$, and thus of the resolution, is 
a product family. For fixed $k$, it follows, 
as in the proof of Theorem \ref{th27900} that by taking
$l \gg 0$, we may assume that $ f_i \in (x,y)^{k+1}$ 
for all $i$. We take $k$ so big, that 
$(x,y)^{k+1} \subset (x,y)^2\left(\frac{\partial f}{\partial x},
  \frac{\partial f}{\partial y}\right)$.  We write $x_1 =x$,
and $x_2 =y$. 
Thus we can  write, formally,  $ \sum_i t^{i-1}f_i
= \sum_i\xi_i(t) \frac{\partial
  f}{\partial x_i}$, with $\xi_i(t) \in (x,y)^2$. Hence
 \[
\frac{\partial F_T}{\partial t}= \sum_i t^{i-1}f_i 
=   \sum_i\xi_i(t) \frac{\partial
  F_T}{\partial x_i}  - \sum_j\sum_i\xi_i(t) t^j\frac{\partial
  f_j}{\partial x_i}. 
\] 
Note  that $\sum_j\sum_i\xi_i(t) t^j\frac{\partial
  f_j}{\partial x_i} \in (x,y)^{k+2}$. 
Iterating this procedure we get
\[
\frac{\partial F_T}{\partial t} = \sum \bar{\xi}_i\frac{\partial F_T}{\partial x_i}.
\]
for some formal $\bar{\xi}_i(t) \in (x,y)^2$. By Artin's Approximation
Theorem,
\cite{art} we may
even assume that the $\xi_i$ are analytic. 
By the characterization of Local Analytic Triviality,
see for example \cite{JP}, Theorem 9.1.7,
it follows that the deformation of $C_T$, given by $F_T= 0$,
is trivial.
By looking at the conductor, it follows that the family $\Sigma_T$ is
trivial. In particular, the induced first order R.C.--deformation
of $(\Sigma,C)$ is trivial. Hence we get a map
\[
 T^1_{X(C,l)} \lra T^1(\Sigma,C),
 \]
  which one easily sees to be linear. 
As $l \gg 0$, it is also surjective, see \cite{djvs2},
 Proof of Theorem 3.3. 
 The kernel one gets, as described in \cite{djvs2},
as follows. Consider   elements
$g ' \in \C\{x,y\}$ such that $g_\eps :=g + \eps g'$ give an R.C--deformation of
$(\Sigma, (g= 0))$, by keeping $\Sigma$ fixed.
As proved in \cite{djvs2}, these are exactly the elements in $I^{ev}$.
Then consider $Y_\eps$ defined by $zf -g_\eps = 0$. 
We get $X_\eps$ as above. Obviously the $g'$ in $\Theta_\Sigma(g)$ 
give trivial deformations $Y_\eps$ of $Y$, and thus trivial
deformations
$X_\eps$ of $X$. 
As soon as $C$ is singular, and the class of $g'$ in
$I^{ev}/(f,\Theta_\Sigma(g))$
is nonzero, this deformation is non-trivial. Then if it where, the
deformation $Y_\eps$ could be extended to a trivial family $Y_T$ given
by $zf -g_T =0$. But as the vanishing order of $g'$ on 
the normalization $\tilde{C}_i$ of at least one branch $C_i$ 
at the point mapping to the singular point of $C$ is
strictly smaller that $c(i) + l(i)$, this then also holds for 
a general $g_t$ for $t$ small. Thus the resolution
graph of the normalization of $zf -g_t$ is different from the
resolution graph of $X(C,l)$. In particular the deformation
is not trivial, contradiction.
Thus we proved the existence of the exact sequence 
(\ref{eq28900}). 

\begin{proof}[Proof of Proposition \ref{Prop27900}]
The proposition follows from the exact sequence
(\ref{eq28900}).
\end{proof}

\begin{remark}
Given the decorated curve $(C,l)$, the kernel 
of the map
\[
\ES_C(\C[\eps]) \lra \ES_{X(C,l)}(\C[\eps])
\]
can be computed explicitly. Namely, one 
searches for equisingular deformations of $C$, which 
deform $X(C,l)$ trivially. 
Given the projection $Y(C,l)$, say given 
by $zf -g = 0$, one can compute all R.C.--deformations
of $Y(C,l)$ which deform $X(C,l)$ trivially.
 Consider generators $u_0=1, u_1,\ldots, u_k$
of $\ko_{X(C,l)}$ as $\ko_{Y(C,l)}$--module. 
In \cite{djvs3}, and \cite{djvs2}, A.9,   for any vector field $\theta$ 
on $\C^3$  an action of $u_i \theta$ on 
$zf-g$ is defined, which, by simultaneously normalizing,
give all infinitesimal trivial deformations of $X(C,l)$. 
These deformations in general, do not keep 
the form $zf-g$ fixed, that is, they are in general
not of this the type $zf_\eps -g_\eps = 0$, 
for some deformation $f_\eps$ and $g_\eps$. 
But we can look at the subspace of deformations
that  do! Thus we get can compute 
all $f'$ so that there exist a $g'$ 
so that $z(f + \eps f') - (g +\eps g')$
that gives an R.C.--deformation of $Y(C,l)$
which gives trivial deformation of $X(C,l)$.
These give equisingular deformations $f + \eps f'$ of
$C$.  
The totality of such $f'$ build an ideal
in $\ES(\C[\eps])$ giving the kernel of the 
map $
\ES_C(\C[\eps]) \lra \ES_{X(C,l)}(\C[\eps])$.

We consider an example. 
For $C$  we take four lines through
the origin, and for $l$ we take the function
assigning two to each branch. The sandwiched
singularity $X(C,l)$ is isomorphic to 
the cone over the rational normal curve of degree $5$.
As this singularity is taut, we get $\es(X(C,l))=0$. 
$Y(C,l)$ can be given by $z(y^4 -x^4) -x^5 = 0$. 
As described in \cite{djvs3}, the generators 
of  $\ko_{X(C,l)}$ as $\ko_{Y(C,l)}$--module,
say $1, u_1, u_2, u_3$, 
correspond to the rows of the matrix:
\[
\left(\begin{array}{cccc} 
zy &0&0& zx + x^2 \\
x & y & 0 & 0\\
0 & x & y & 0  \\
0 & 0 & x & y
\end{array}
\right).
\]
Note that the lower three rows gives a resolution
of the conductor of $C$,
and that $zf-g = z(y^4 -x^4) -x^5 = 0$ 
is the determinant of this matrix. 
The columns give linear equations, so we have four of them:
\[
\begin{array}{ll}
L_1:& zy + xu_1 = 0\\
L_2: &yu_1 + xu_2 = 0\\
L_3: &yu_2 + xu_3 =0\\
L_4: & (zx+x^2) + yu_3 =0.
\end{array}
\]
From these, the quadratic equations can be calculated:
\[
\begin{array}{l}
u_1^2 = zu_2\\
u_2^2 = z(z+x)\\
u_3^2 = (z+x)u_2\\
u_1u_2 = zu_3\\
u_1u_3 = z(z+x)\\
u_2u_3=(z+x)u_1.
\end{array}
\]
 Thus we get a total of ten equations, describing 
the cone over the rational normal curve of degree $5$.
The actions of $u_i \theta$ on these equations are
totally determined by their values on the linear
ones, but one needs the quadratic equations to compute
them. We do not want to do the whole calculation in detail,
but only look at the action of $u_3 \tfrac{\partial}{\partial x}
-z \tfrac{\partial}{\partial y}$. We get the values
\[
\begin{array}{l}
L_1 \mapsto u_1u_3 -z^2 = xz\\
L_2 \mapsto u_2u_3 -zu_1  = xu_1\\
L_3 \mapsto u_2u_3 -  u_3u_2 =0\\
L_4 \mapsto zu_2 +2xu_2 - u_3^2 = xu_2.
\end{array}
\]
Thus we get the following infinitesimal deformation of 
the matrix:
\[
\left(\begin{array}{cccc} 
z(y+ \eps x)  &0&0& zx + x^2 \\
x & y + \eps x & 0 & 0\\
0 & x & y & \eps x  \\
0 & 0 & x & y
\end{array}
\right)
\]
which gives the R.C.--deformation:
\[
z (y^4 -x^4 + \eps x^2y^2) - x^4.
\]
Thus we see that  the non--trivial equisingular deformation
$y^4 -x^4 + \eps x^2y^2$
of $C$
maps to the trivial infinitesimal deformation
of the cone over the rational normal curve of degree $5$. 
\end{remark}

\section{The Dimension of the Artin Component}

We quote the following result.

\begin{Theorem}[Wahl, \cite{wahl1}, Propositions 2.2 and 2.5]\label{wa}
Let $X$ be a rational surface singularity,
$p: \tilde{X} \lra X$ be the minimal resolution, 
 $E$  be the exceptional divisor 
of $p: \tilde{X} \lra X$, and
$E_1, \ldots, E_s$ be the irreducible components of $E$. Let $b_i =
E_i^2$. 
Then
\[
 h^1(\tilde{X}, \Theta_{\tilde{X}}) = -\sum_{i=1}^s(b_i+1) +
 h^1(\tilde{X}, \Theta(\log E)).
\]
The vector space  $H^1(\tilde{X}, \Theta(\log E))$ is the Zariski 
tangent space to the base space of the equisingular deformations
of $X$. Its dimension we denote by $\es(X)$ for short.
\end{Theorem}

The space $H^1( \tilde{X}, \Theta_{\tilde{X}})$ classifies
the infinitesimal deformations of the resolution $\tilde{X}$.
Its dimension is the dimension of  the Artin component. We will give
a different characterization
of this dimension for sandwiched singularities
$X(C,l)$, for the case $l \gg 0$. 
Note moreover that we can calculate $I^{ev}/(f,\Theta_\Sigma(g))$
on the normalization, as the conductor $I$
is contained in $I^{ev}$. 
Thus we get that
\[
\dim_\C\left( I^{ev}/(f,\Theta_\Sigma(g)\right) = \sum_{i \in T} \left(l(i)
- m(i)\right).
\] 
Hence from the exact sequence (\ref{eq28900}) it follows that
\begin{Theorem}
Let $(C,l)$ be a decorated curve,
and suppose that $l \gg 0$. Let $\Sigma$ be the conductor of $C$. Then
\[
\dim_C(T^1_{X(C,l)}) = \sum_{i \in T} \left(l(i)
- m(i)\right) + \dim_\C(T^1(\Sigma,C)).
\]
\end{Theorem}

We want to stress that in general for small $l$ the statement
of the theorem is false.\\ 

We want to understand the deformations on the Artin component
of $X(C,l)$. This has been described in \cite{djvs2}, 4.13. 
One can decide whether a
  deformation occurs on the Artin component by looking at 
the corresponding R.C.--deformation of $(\Sigma,C)$. 
A general such one--parameter  R.C.--deformation 
of $(\Sigma,C)$ is one, which has on a general fiber $q$
singular points, where $q$ is the number of infinitely near points
of $C$. For each infinitely near point $P$ of $C$ of multiplicity
$m_P$
we have on a general fiber of the deformation a singularity
consisting of $m_P$ smooth branches intersecting transversely
(which we call ordinary $m_p$--tuple point). 
Let $A$ be the closure of the stratum of the base space
of R.C.--deformations of $(\Sigma,C)$ where the above mentioned
deformations occur.  
All elements of $I^{ev}$ are unobstructed 
against these deformations.
This follows immediately from the picture method, see \cite{djvs2}.
As the Artin component is smooth, it follows that $A$ is smooth. 
We get:

\begin{Proposition}\label{prop28900}
Let $(C,l)$ be a decorated curve with $l \gg 0$. Let $A$ be the
stratum described above. Then the dimension of the Artin component 
of $X(C,l)$ is equal to
\[
\dim(A) + \sum_{i \in T}\left( l(i)-m(i)\right).
\]
\end{Proposition}
We now compute $\dim(A)$. The result is:
\begin{Proposition}\label{prop28900A}
Let $C$ be a plane curve singularity, and $A$ be the stratum described
above. Then
\[
\dim(A) = \tau(C) -\sum_{P \in \IN(C)}\tfrac{1}{2}(m_P^2 + m_p -4).
\]
\end{Proposition}

\begin{proof}
We have an R.C. deformation of $(\Sigma, C)$ over $A$.
In particular we get a one--parameter deformation of $C$ over $A$. 
Let $B$ be the base space of a semi-universal deformation 
of $C$. It is a smooth space of dimension $\tau(C)$.
By semi-universality we get a map of smooth spaces $A \lra B$.  
This map in general is not an immersion. Indeed, 
Buchweitz \cite{Bu} showed that the kernel of the map $T^1(\Sigma,C) \lra T^1_C$
is equal to $m(C) -r(C)$, the multiplicity minus the number of
branches.
However, for a general point of $A$ we have that for all 
singularities $m_P(C) = r_P(C)$,
so that the map $A \lra C$ is an immersion at a general point of $A$.
Hence the image of  $A$ at a general point is smooth of dimension
$\dim(A)$. Thus, by openness of versality, it remains to 
compute for each infinitely near point $P$, the codimension
of the stratum $A$ of the ordinary $m_P$--tuple point in 
the base space of the ordinary $m_P$--tuple point. 
This  codimension  is $\frac{1}{2}(m_P^2 + m_p -4)$,
being two less than the number of monomials of degree smaller than $m_p$.
\end{proof}

As an application, we can give the dimension of the equisingular 
stratum of a plane curve singularities.
This formula is  equivalent to formulas
given by
Wall \cite{wall} and Mattei, \cite{mattei}.
 A special
case also has been considered in Piene and Kleiman, \cite{piene}.
\begin{Theorem}\label{th28900} The following formula holds for all plane curve singularities $C$:
\[
\es(C) = \tau(C)  + \eb(C) +\sum_{i=1}^k (b_i+1)-\sum_{P \in \IN(C)}
\tfrac{1}{2}(m_P^2 + m_p -4).
\]
where
\begin{enumerate}
\item $\tau(C)$ is the Tjurina number.
\item  $\IN(C)$ is the set of infinitely near points of $C$ (including
the singularity itself); $m_P$ the multiplicity of the 
infinitely near point $P \in \IN (C)$. 
We only consider infinitely near points of $C$
which are singular, that is, $m_P >1$,
\item  the  $b_i$ are  the self-intersection
numbers of the exceptional curves  $E_i$ on a minimal good embedded resolution of $C$.
\end{enumerate}
\end{Theorem}
\begin{proof}
Given $C$, take a decorated curve $(C,l)$ with $l \gg 0$. 
Let $E_1, \ldots, E_s$ be the exceptional curves in 
the minimal resolution of $X(C,l)$.
Let $b_i = -E_i^2$.  By Wahl's result \ref{wa},
the dimension of the Artin component of $X(C,l)$ is
equal to 
\[
\es(X(C,l)) -  \sum_{i=1}^s(b_i+1).
\]
On the other hand, combining \ref{prop28900}
and \ref{prop28900A} we get that this dimension is also equal
to 
\[
\sum_{i \in T} \left(l(i) - m(i)\right) + \tau(C) - \sum_{P \in \IN(C)} \tfrac{1}{2}(m_P^2 + m_p -4).
\]
Using \ref{Prop27900} we thus get
\[
\es(C) = \tau(C) + \sum_{i=1}^s(b_i+1) + \sum_{i \in T}\left(l(i) -m(i)\right)
-\sum_{P\in \IN(C)}\tfrac{1}{2}(m_P^2 + m_p -4).
\]
It remains to show that $\sum_{i=1}^s(b_i+1) + \sum_{i \in T}\left(l(i)
-m(i)\right) = \eb(C) +  \sum_{i=1}^k(b_i+1)$.
This is easy, as we know that we get the minimal resolution of
$X(C,l)$
out of the minimal good resolution of $C$ by doing an 
extra $l(i) -M(i)$ blowing-ups, for each $i \in T$.
We
 thereby introduce a chain of $(-2)$--curves of length $l(i)
-M(i)-1$,
and decrease the selfintersection of the exceptional
curve on the minimal resolution of $C$ which $C_i$ intersects
by one. 
From $\eb(C) = \sum_{i \in T}M(i) -m(i)$,
see \ref{propa}, the result
follows.
\end{proof}

\begin{remark}
It is not so difficult to see that for a decorated
curve $(C,l)$ with $l \gg 0$ we have $m(X(C,l)) = m(C) +1$.
Note that we in fact computed the codimension of the
Artin component in  $T^1_{X(C,l)}$. With a little combinatorics
one then proves that the invariant $c(X(C,l))$, introduced
by Christophersen and Gustavsen \cite{CG}, is equal to  $0$. 
This can also be proven directly, using the results
of \cite{CG} and  \cite{G}, chapter two.
In particular, for $l \gg 0$ one gets a formula for the 
dimension of $T^2_{X(C,l)}$ by using the results
of Christophersen and Gustavsen. By using  semicontinuity
of $T^2$, this then also holds for all $X(C,l)$
with $m(X(C,l)) = m(C) +1$. 
Also, in case $l \gg0$, one can construct, using the picture method, see
\cite{djvs2},  a deformation
of $X(C,l)$, where on a general fiber we have a cone over the rational
normal curve of degree $m_P$ for each infinitely near point $P$ of
$C$.
Using a standard argument, as used in \cite{djvs1}, we
then get that the obstruction map is surjective, that is,
the minimal number of equations for describing the base space
of a semi-universal deformation of $X(C,l)$ is equal
to the dimension of $T^2_{X(C,l)}$ if $l\gg 0$. 
\end{remark}

\section{Computation of the Equisingularity Ideal}

  Let $C$ be a plane curve singularity, and let $s: \Spec(\C) \lra C$ be 
the singular point of $C$. We consider the functor of Artin
rings $E = E_C$ of deformations of $C$ which are equimultiple
along a section, see Wahl \cite{wahl} 1.3. That is,
for an Artinian ring $A$,  an element
of $E(A)$ is a pair $C_A \lra \Spec(A)$,
of deformations of $C$ together with a deformation  $\bar{s} : \Spec(A) \lra C_A$
of $s$ such that $C_A$ is equimultiple along $\bar{s}$. 

 Next Wahl considers the functor $E^{(2)}$.
An element of $E^{(2)}(A)$ consists of 
an element  $C_A \lra \Spec(A),$ $\bar{s} : \Spec(A) \lra C_A$
of $E(A)$,
 together with sections 
$\bar{s}_i : \Spec(A) \lra B_{\bar{s}}$ for $i = 1,\ldots,t$.
Here $B_{\bar{s}}$ is the blowing up along the section $\bar{s}$, and 
$\bar{s}_i : \Spec(A) \lra B_{\bar{s}}$ should
induce equimultiple deformations of all singular points
of the {\em reduced total transform}  of $C$ under
blowing up the origin in $\C^2$.
 Inductively one  defines the functor $E^{(n)}$,
see Wahl \cite{wahl}, 2.7.

We will also  consider a slightly different functor
$E_{(2)}$. Elements of  $E_{(2)}(A)$ are like elements
of $E^{(2)}(A)$ except that for non-nodes
of the reduced total transform  the $\bar{s}_i : \Spec(A) \lra B_{\bar{s}}$ 
induce equimultiple deformations of the intersection points
of the {\em strict  transform} $\bar{C}$ of $C$ and the exceptional
curve 
 under
blowing up the origin in $\C^2$. Inductively we define
the  functor   $E_{(n)}$:
\begin{defn}
Let $C$ be a plane curve singularity, and let $\bar{C}$ for
$i=1, \ldots t$ be the connected components of the strict
transform of the first blowing up of $C$ (corresponding
to the tangential components of $C$). 
Suppose that for $2 \leq j \leq n-1$ the functor $E_{(j)}$ 
has been defined, together with natural maps $E_{(j)} \lra E_{(j-1)}$.
Then an element of $E_{(n)}(A)$ by definition is
a $(\alpha; \bar{u}_1, \ldots, \bar{u}_m)$ with
 \begin{enumerate}
\item 
$\alpha \in E_{(n-1)}(A)$. We let $B_\alpha$ be the $A$--space
obtained by blowing up successively along the $A$--sections of
$\alpha$. We let $\alpha^* = (C_A, \bar{s}; \bar{s}_1,
\ldots, \bar{s}_t)$ be the image of $\alpha$  in $E_{(2)}$.

For each $j=1, \ldots, t$, we therefore have
ordered collections
$\alpha_j$ of $A$--sections in $\alpha$ such
that $(\bar{C}_j,\alpha_j) \in  E_{n-1\bar{C}_j}(A)$.
\item
$\bar{u}_i$ is an $A$--section of $B_\alpha$.
\item Every $\bar{u}_i $
 lies over a section $\bar{s}_j$ of $B_{\bar{s}}$ in $\alpha^*$;
 letting $\bar{u}_{j,1}, \ldots ,\bar{u}_{j,q} $ denote all such
sections
we have for all $j=1, \ldots, t$ that
\[
(\bar{C}_j , \alpha_j; \bar{u}_{j,1}, \ldots ,\bar{u}_{j,q}) 
\in  E_{n-1\bar{C}_j}(A).
\]
\end{enumerate} 

\end{defn}
 Exactly as in Wahl
 2.5 and 3.2 one shows that $E_{(n)}$ has a very good deformation
theory, and that for $N \gg 0$ the natural maps $E_{(N+1)} \lra E_{(N)}$
is bijective. It is however {\em false } in general that for
such a big $N$
the natural morphism  from $E_{(N)}$ to the functor of deformations
of $C$ is injective, contrary to the functor $E^{(N)}$.

Wahl defines the functor of equisingular deformation
$\ES$ to be $E^{(N)} $ for $N \gg 0$. It is a subfunctor
of the deformation functor of $C$. Similarly we define:
\begin{defn}
For a plane curve singularity we define the 
 simultaneous resolution functor $\SR = \SR_C$ to be
$E_{(N)}$ for $N \gg 0$.
 \end{defn}
 The functor $\SR$  is in general not a subfunctor
of the deformation functor of $C$. 
We now calculate the  tangent space $\SR(\C[\eps])$ 
of $\SR$ inductively.
For this,  consider 
 the  blowing up
$
(Y,E) \stackrel{p}{\lra} ( \C^2, 0).
$
The strict transform of $C$ under $p$ split
 up into say $t$ connected components, say $\bar{C}_1, \ldots, \bar{C}_t$.
We define $C_i = p(\bar{C}_i)$. We suppose that $C$ is given by $f =0$ for 
a square free $f$, and we let $f = f_1\cdots f_t$,
where $C_i$ is given by $f_i=0$.
We let $p_i$ be the intersection point of $\bar{C}_i$ with
$E$.

We chose coordinates $x,y$ on  $\C^2$, so we get homogeneous
coordinates $(u:v)$ on the exceptional projective line $E$.
The blowing up is given by the equation $uy =xv$.  For simplicity,
 we  will
 assume that none of the $p_i$ is the point with
homogeneous coordinates $(0:1)$. 
The homogeneous coordinates of $p_i$ we denote
 by $(1:a_i)$. Thus we only have to look
at the chart $u =1$ in the blowing up. In these coordinates,
the blowing up is given  by $y=xv$, and thus the strict transform
 $\bar{C}$ is given by  $\bar{f}(x,v)= 0$, where
we have the equality $f(x,y) = f(x,xv) = x^m \bar{f}(x,v)$.
 Here $m$ is the multiplicity of $C$. Similar for $f_i$ and
$\bar{f}_i$.
 
 By induction we may suppose known the tangent
 space $\SR_{C_i,p_i}(\C[\eps])$ for 
 $i =1, \ldots ,  t$. 
  For each element in $\SR_{C_i,p_i}(\C[\eps])$
 we get, in  particular, an
 infinitesimal (equimultiple) deformation of $(C_i,p_i)$
 given, say, by   $f_i +\eps g$.
 Thus we have a natural map $\SR_{C_i,p_i}(\C[\eps]) \lra
\C\{x,v-a_i\}$,
 whose image we call $J_i$. 
 If $m_i$ is the multiplicity of $C_i$, we may assume that
the curve $C_i$  is given by a Weierstrass polynomial
of degree $m_i$ in $v$.  Applying the Weierstra{\ss}
Division Theorem,  we may assume that the degree of 
$g$ in $v$ is smaller  than $m_i$.
  
Now  let $((f_i + \eps \bar{h}_i), \alpha_i)$ be 
 a   element of $\SR_{C_i,p_i}(\C[\eps])$. Thus
$\alpha_i$  is an ordered system of infinitesimal sections, and we assume that 
  $\bar{h}_i(x,v-a_i)$
  has 
  degree smaller than $m_i$ in $v$.  We consider
  $x^{m_i}\bar{h}_i(x,v-a_i)$, where $m_i$ is the multiplicity of $C_i$.
   Using the equation $y=xv$, we can eliminate $v$, and get 
  an  element  $h_i(x,y) \in \C\{x,y\}.$
 We  consider the product $h_i'=h_i(x,y)\prod_{j\neq i}f_i$,
 and  the sum $h'=\sum_{i=1}^t h_i'$. 
 Together  with the trivial section $\bar{s}$
 we get the  following element
\[
 \left((f + \eps h'), \bar{s},\alpha_1, \ldots, \alpha_t\right).
\] 
  of $\SR(\C[\eps])$.
We thus get a map 
$\psi: \prod_i\SR_{C_i}(\C[\eps]) \lra  \SR_{C}(\C[\eps]).$
 The following is obvious. 
 
\begin{Lemma}
The image of $\psi$ consist of  all elements of $\SR(\C[\eps])$ whose
first infinitesimal section is trivial.
\end{Lemma}

 Furthermore, we have the trivial deformations of $C$,
which lie in the Jacobian ideal of $f$. They certainly
 resolve simultaneously. Those for which the 
first section is not trivial, 
have as image in $E$ the infinitesimal deformations
$\bigl((f + a\eps\frac{\partial f}{\partial x} +
 b\eps\frac{\partial f}{\partial y}), (x+a\eps,x +b\eps)\bigr)$
for $(a,b)\in \C^2\{(0,0)\}$. These lead to 
two elements $\alpha_x$ and $\alpha_y$ in  $\SR(\C[\eps])$,
corresponding to the pairs $(a,b) = (1,0)$ and $(0,1)$. 
 \begin{Proposition}
 For a plane curve singularity, $\SR(\C[\eps])$ is
generated by the image of $\psi$, $\alpha_x$
and $\alpha_y$. 
\end{Proposition}

\begin{proof}
Let an element $\alpha$ in $\SR(\C[\eps])$ whose first
 section is given by the ideal $(x-a\eps, y-b\eps)$ 
for $(a,b) \in \C^2$. Then $\alpha +a\alpha_x +b \alpha_y$
has trivial first section, hence is contained in the image
of $\psi$ by the lemma. 
\end{proof}

 We can easily calculate generators of the image of $\psi$,
and thus of $\SR(\C[\eps])$ as follows. 
Namely, take generators of $\SR_{C_i,p_i}(\C[\eps])$
  such that their images in $J_i$ (say under
the map $\beta$) generate
 $J_i$ as $\C\{x\}$--module. 
That is, in each degree $p <m_i$ we calculate all 
elements of $J_i$ of degree $p$ in $v$.
As the polynomials in $v$ of degree smaller than
$m_i$ is a finitely generated $\C\{x\}$--module, 
and $J_i$ is a submodule, we can calculate
 such generators by computing a standard basis.
The inverse image under $\beta$ of such a generator
is  a finitely generated $\C\{x\}$--module. 
By induction one only has to show this for an infinitesimal equimultiple
deformation. But the different sections form a 
two dimensional vector space, so 
certainly finitely generated $\C\{x\}$--module. 

So take finitely many of those generators
of $\SR_{C_i,p_i}(\C[\eps])$ as $\C\{x\}$--module. Their images under $\psi$
then generate the image of $\psi$.\\

Now we turn our attention to the calculation of $\ES(\C[\eps])$. We 
consider the minimal good resolution $(Z,F)$ of $C$,
 and we let $F_1, \ldots. F_p$ be the irreducible components of $F$.

\begin{Theorem}\label{th291200}
There is a natural map of vector spaces
\[
\vfi: \SR(\C[\eps]) \lra \oplus_{i=1}^p H^1(F_i, N_{F_i/Z})
\]
whose kernel is $\ES(\C[\eps])$: the tangent space of the equisingular
deformations of  $C$. 
\end{Theorem} 

\begin{proof}
 Clearly we have  
$\ES(\C[\eps]) \subset \SR(\C[\eps])$.
 Conversely, suppose $\alpha$ is an element of $\SR(\C[\eps])$.
Associated to it is a space $B_\alpha$, which
is obtained by blowing up $\C^2 \times \Spec(\C[\eps])$
successively along the sections described by $\alpha$.
Hence we get an infinitesimal deformation of some resolution
of the curve $C$, in particular a deformation of the minimal
resolution  of $Z$, that is
an element of $H^1(Z,\Theta_Z)$. We look at the induced deformation
of a small neighborhood of each $F_i$. More precisely,
 from the natural composition
 $\Theta_Z \lra \oplus \Theta_Z \tensor \ko_{F_i} \lra \oplus N_{F_i/Z}$ we get
a map $H^1(Z,\Theta_Z) \lra \oplus_{i=1}^p H^1(F_i, N_{F_i/Z})$. Thus
gives the map $\vfi$. The kernel consists of elements of $\SR(\C[\eps])$
  which induce
trivial deformations of all small neighborhoods of each $F_i$. 
This is equivalent to saying that we get an
infinitesimal equisingular deformation of the minimal resolution $Z$,
that is, an element of $\ES_Z(\C[\eps])$. By
Wahl, \cite{wahl}, 5.7  we have a natural smooth
functor $\ES_C \lra \ES_Z$, whose kernel on tangent spaces  is given by the
tangent space of the functor 
$\ES'$, which is the functor of all equisingular deformations 
 of $C$ for which all sections can be trivialized, 
see \cite{wahl} 5.4. This proves the Theorem. 
\end{proof}

\begin{remark}
Note that the   $\ES'(\C[\eps])$, the tangent space
of all equisingular deformations of $C$ whose sections
can be trivialized, certainly is in the kernel of $\vfi$.
In this case, the resolution $Z$ is  deformed trivially. 
This $\ES'(\C[\eps])$ is  easy to compute,
see Wahl \cite{wahl} 6.3. We thus get an induced map
\[
\vfi: \SR(\C[\eps])/ \ES'(\C[\eps]) \lra \oplus_{i=1}^p H^1(F_i, N_{F_i/Z})
\]
\end{remark}

Given an element $\alpha$ in $\SR(\C[\eps])$,
one can write down the space $B_\alpha$,
which then gives the wanted
 infinitesimal deformation of $Z$. 
We have to see what the image in each $H^1(E,N_{E/Z})$ is
for each exceptional curve $E$. 
We first recall how to describe  $H^1(E,N_{E/Z})$ by
means of Cech cohomology. 
In our case, the embedding $E$ in $Z$ will
always be described in the following way. On $E$ we have 
 homogeneous coordinates $(u:v)$.
The space $Z$ is 
in a neighborhood of $E$, given by two charts.

\begin{itemize}
\item In the first chart $U_1$ we have coordinates $v$ and $x$.
\item In the second chart $U_2$ we have coordinates $y$ and  $u$.
\item The transition functions on the intersection $U_1\cap U_2$ are
\[ 
u = v^{-1}, \quad y = \alpha(v) x
\]
where $\alpha(v)$ is a polynomial in $v$ of degree $k$ (with zero constant
coefficient). 
\end{itemize}
The curve $E$ is given by the zero set of  $I:=(x, y)$.
In fact, if $Z$ is obtained from the blowing up of $\C^2$, by
blowing up points $(1:a_1), \ldots, (1:a_{k-1})$ (not necessarily distinct),
 and no other points on $E$, 
then the polynomial $\alpha(v)$ is equal
to $v \prod_{i=1}^{k-1}(v-a_i)$. 

The Cech cohomology group $H^1(E, N_{E/Z})$ is easily described. 
 On the intersection $U_1 \cap U_2$, it is given
by an element of $\Hom(I, \ko/I)$ which sends
 $y$ to $h(v)$, where $h(v)$ is a polynomial of degree smaller
 than $k$ with no constant coefficient. (Note that $v= u^{-1}$ on $U_1 \cap U_2$.)
If it sends $y$ to an arbitrary holomorphic function $h(v)$, then the class
in  $H^1(E, N_{E/Z})$ is obtained by the remainder of $h(v)$ through
 $v \prod_{i=1}^{k-1}(v-a_i)$, and forgetting the constant term. 
The corresponding infinitesimal deformation of a small neighborhood
of $E$ in $Z$ is given by the transition functions
\[ 
u = v^{-1}, \quad y = \alpha(v) x + \eps h(v).
\]
 This is the Kodaira-Spencer description of infinitesimal deformations. 
 So in this way we can make the map $\vfi$ explicit,
and thus calculate the equisingularity ideal.\\

  In practice it is probably better to proceed with induction:
suppose known the $\ES_{\bar{C}_i}(\C[\eps])$
for $i=1, \ldots, t$, where the $\bar{C}_i$ are the connected components
    of the strict transform of $C$ under blowing up the 
origin in $\C^2$. (Note that these can all be identified with ideals,
that is, the sections are uniquely determined by the deformation.)
Let $E$ be the exceptional divisor, and suppose that
on 
the minimal good resolution $(Z,F)$ the curve  $E$ has
selfintersection $-k$. That means, that
we have to blow up $k-1$ times in the 
exceptional curve of the blowing up up 
 $\C^2$ in the process to arrive at the minimal
 good resolution $(Z,F)$. 
We consider the submodule $\AES :=\psi (\prod_{i=1}^t \left(
  \ES_{\bar{C}_i}(\C[\eps]) \right)
\subset \SR(\C[\eps])$.

\begin{Lemma}
The map $\AES \lra  \oplus_{F_i \neq E} H^1(F_i,N_{F_i/Z})$
induced by $\vfi$ is the zero map. 
\end{Lemma}

\begin{proof}
 After blowing up in the trivial section, we get,
by {\em construction} of $\AES$  equisingular deformations 
 of all connected components of the strict transforms of $C$,
which then, by Theorem 
\ref{th291200}, map to the zero element of   $\oplus_{F_i \neq E} H^1(F_i,N_{F_i/Z})$.
\end{proof}

We therefore get an induced map
\[
\vfi: \AES  \lra  H^1(E,N_{E/Z}),
\] 
 whose kernel is equal to $\ES(\C[\eps])$. 
 Obviously $\ko_C/\ES'(\C[\eps])$ is a finite dimensional
 vector space. Thus we get an induced map
 of finite dimensional vector spaces:
\[
\vfi: \AES /\ES'(\C[\eps])  \lra  H^1(E,N_{E/Z}).
\]

\begin{example}
We take the curve singularity given by
\[
f = (y^4 + x^5)^2 +x^{11} = 0. 
\]
 The Tjurina number $\tau(f)$  one computes 
 with {\sc Singular} \cite{sing} to be $55$. 
The dimension of the equisingular stratum one computes to be 
$6$.
   The total transform of the blowing up is given by
$x^8\bigl((v^4 +x)^2 x^3=0$.
 We see that the strict transform given by $\bar{f}=(v^4 +x)^2 +x^3=0$ is
   an $A_{11}$ singularity, whose equisingularity
 ideal is trivial, that is, generated by the partial derivatives
 of a defining equation. Thus it
 is generated as $\C\{x,v\}/(\bar{f})$--module by 
$\frac{\partial \bar{f}}{\partial x}$ and $\frac{\partial
  \bar{f}}{\partial v}$.
 We need however generators as $\C\{x\}$--module. $\C\{x,v\}/(\bar{f})$
 is a free $\C\{x\}$--module of rank $8$, generated
 by $1, v, \ldots, v^7$. Thus we have the following
  $\C\{x\}$--generators:
 \begin{eqnarray*}
 \frac{\partial \bar{f}}{\partial x}, v\frac{\partial \bar{f}}{\partial
   x},
 \ldots, v^7\frac{\partial \bar{f}}{\partial
   x},\\
 \frac{\partial \bar{f}}{\partial v}, v\frac{\partial \bar{f}}{\partial
    v},
  \ldots, v^7\frac{\partial \bar{f}}{\partial v}.
\end{eqnarray*}
 We need classes of these elements in $\C\{x,v\}/(\bar{f})$ 
   of degree smaller than eight in $v$. For these we can take
\[
2(v^4 +x) + 3x^2,  \ldots, 2v^3(v^4 +x)+3v^3x^2
\]
\[
 -2x(v^4 + x +x^2) + 3v^4x^2, \ldots,  v^3(-2x(v^4 + x +x^2) +
3v^4x^2),
\]
\[
v^3(v^4 + x), xv^4 + x^2 + x^3, \ldots, v^3(xv^4 + x^2 + x^3),
\]
\[-x^2(v^4 + x + x^2) + v^4x^3 \ldots,
-x^2v^2(v^4 + x + x^2) + v^6x^3.
\]
We first look at
$\frac{\partial \bar{f}}{\partial x} = 2(v^4 +x) + 3x^2$.
The corresponding section is given by $(x +\eps,v)$.
Thus we blow up the 
trivial infinitesimal deformation of the first blowing up given by the transition functions
\[
 u = v^{-1}, \quad y = v x
\] 
 in the section $(x+\eps,v)$.  In the appropriate chart
this blowing up is given by the equation $x+ \eps = vx_1$. 
Hence after this blowing up, a neighborhood of $E$ is
given by the transition functions
\[
 u = v^{-1}, \quad y = v^2x_1 - \eps v.
\] 
On the first blowing up, the equisingular deformation
of the $A_{11}$ singularity was given
by $\bar{f}( x +\eps,v) = 0$. If the second blowing up
of the undeformed curve is given by $\bar{\bar{f}}(x_1,v)=0$, we see that 
for the deformed one we have
$\bar{\bar{f}}(x_1,v)=0$. Hence we
can continue blowing up in trivial sections, leading
to the deformation of $Z$ given in a neighborhood of $E$ by the transition functions
\[
 u = v^{-1}, \quad y = v^5x_4 - \eps v.
\] 
This gives a nontrivial element in $H^1(E, N_{E/Z})$. 

Now look at $v \frac{\partial \bar{f}}{\partial x}$.
The corresponding infinitesimal deformation of the strict transform of $C$  is given by
 $\bar{f}( x +\eps v,v) = 0$. We blow up in the trivial section $(x,v)$
and in the appropriate chart we have $x =vx_1$, so we get 
\[
 u = v^{-1}, \quad y = v^2x_1.
\] 
The induced deformation on this blowing up of $C$ is 
described by  $\bar{f}( v x_1+\eps v,v)= v^2\bar{\bar{f}}(x_1 + \eps,v)$.
Thus we can apply the previous discussion and blow up in $x_1 + \eps$ leading
to the third blowing up given by $x_1 + \eps = vx_2$ getting the transformation
functions
\[
 u = v^{-1}, \quad y = v^3x_2 -\eps v^2.
\]
We can now, as before, continue with trivial sections
leading to an infinitesimal deformation of $Z$ given in a neighborhood of $E$ by the transition functions
\[
 u = v^{-1}, \quad y = v^5x_4 - \eps v^2.
\] 
This  again gives a nontrivial element in $H^1(E, N_{E/Z})$. 
In general one see that the deformation given by
\[
\bar{f}(x +\eps g(v),v)
\]
for $g$ of degree less than $4$
leads to an infinitesimal deformation of the minimal resolution of $Z$ which
in a neighborhood of $E$ is described by the transition functions
\[
 u = v^{-1}, \quad y = \alpha(v)w -\eps vg(v).
\]
For terms of degree greater or 
equal to $4$, we get trivial deformations along $E$, thereby
inducing infinitesimal equisingular deformations of $C$.

Thus for example $v^4 \frac{\partial \bar{f}}{\partial x}=
-2xv^4 -2x^2-2x^3 +3x^2v^4$ gives the element
 $x^8(-2xv^4 -2x^2-2x^3 +3x^2v^4) = -2x^5y^4 -2x^{10}-2x^{11}
 +3x^6y^4$.
This gives an equisingular deformation. For 
 \[
v^5 \frac{\partial \bar{f}}{\partial x}, \quad
v^6 \frac{\partial \bar{f}}{\partial x},\quad 
v^7 \frac{\partial \bar{f}}{\partial x},
\]
we get
\[
-2x^4y^5 -2x^{9}y-2x^{10}y
 +3x^5y^5, -2x^3y^6 -2x^{8}y^2-2x^{9}y^2
 + 3x^4y^6, 
\]
\[-2x^2y^7 -2x^{7}y^3-2x^{8}y^3
 +3x^3y^7.
\]
Similarly one sees that the deformations
corresponding to  $\frac{\partial \bar{f}}{\partial v},
\ldots, v^7 \frac{\partial \bar{f}}{\partial v}$ all lead
to equisingular deformations of $C$.
These we see to map under $\psi$ to
\[
x\frac{\partial f}{\partial y}, \;
x^5y^4 + x^{10} + x^{11},\:
x^4y^5 + x^{9}y + x^{10}y,\;
x^3y^6 + x^{8}y^2 + x^{9}y^2,
\]
\[
x^2y^7 + x^{7}y^3 + x^{8}y^3,\;
-x^6y^4 - x^{11}  - x^{12} + x^7y^4,\;
 -x^5y^5 - x^{10}y  - x^{11}y + x^6y^5,
\]
\[ 
 -x^4y^6 - x^{9}y^2  - x^{10}y^2 + x^5y^6. 
\]
With {\sc singular} we can now compute 
a standard basis for the equisingularity ideal 
of $C$ (as $\C\{x,y\}$--module). The result is
\[
y^7+x^5y^3, \:
10x^4y^4+10x^9+11x^{10}, \:
x^3y^6+x^8y^2+x^9y^2, \:
x^{11}, \:
x^{10}y, \:
x^9y^2, \:
x^8y^3,
\]
and that indeed  the equisingular stratum has dimension $6$.
The equisingular deformations of $C$ 
all of whose sections are trivial one computes 
to be the ideal
\[
(x^6,x^5y,x^3y^2,x^2y^3,xy^4,y^5)^2 \mod f.
\] 
Together with the trivial deformations these give
only a five dimensional stratum.  The ``extra''
infinitesimal  equisingular deformation, is given by
\[
(y^4 + x^5 + \eps x^3y^2)^2 + x^{11}.
\]
\end{example}

\begin{remark}
Luengo informed me that he has  
a different method for computing the equisingularity 
ideal. Furthermore, Campillo and Greuel 
seem to have a method for computing the equisingularity
ideal 
by using the Hamburger 
Noether expansion. 
A more detailed study of the map $\vfi$
leads to an another proof of Theorem \ref{th28900}.
\end{remark}

\end{document}